\documentclass[preprint,12pt]{elsarticle}

\usepackage[utf8]{inputenc}
\usepackage[francais]{babel}
\usepackage[T1]{fontenc}
\usepackage{amsmath}
\usepackage{amsfonts}
\usepackage{amssymb}
\usepackage{graphicx}
\usepackage{multicol}
\usepackage{multirow}
\usepackage{booktabs}
\usepackage{colortbl}
\usepackage{tabularx}
\usepackage{multirow}
\usepackage{threeparttable}
\usepackage{etoolbox}
\appto\TPTnoteSettings{\footnotesize}
\addto\captionsfrench{}

\usepackage{amssymb}
\usepackage{amsmath}
\newtheorem{thm}{Theorem}

\newtheorem{cor}[thm]{Corollary}
\newdefinition{rmk}{Remark}
\newdefinition{defn}{Definition}
\newproof{pf}{Proof}
\newproof{pr}{Proof}
\newproof{pot}{Proof of Theorem \ref{thm2}}
\newtheorem{alg}[thm]{Algorithm}
\newcommand{\norm}[1]{\left\lVert#1\right\rVert}

\journal{}

\begin{document}

\begin{frontmatter}



%
 \title{
 A treatment of  breakdowns and near breakdowns in a reduction of a matrix to upper $J$-Hessenberg form and related topics	\\
}
 \author[AS]{Ahmed Salam\corref{cor1}}
 \ead{ahmed.salam@univ-littoral.fr}
 \author[AS,HB]{Haithem Ben Kahla}
 \ead{benkahla@univ-littoral.fr} 
 \cortext[cor1]{corresponding author}

 \address[AS]{University Lille Nord de France, ULCO, LMPA. BP 699, F62228, Calais, Cedx, France.}
 \address[HB]{University of Tunis El Manar, ENIT-LAMSIN, BP 37, 1002, Tunis, Tunisia.}



\begin{abstract}
The reduction of a matrix to an upper $J$-Hessenberg form is a crucial step in the $SR$-algorithm (which is a $QR$-like algorithm), structure-preserving, for computing eigenvalues and vectors, for a class of structured matrices. 
 This reduction may be handled via the algorithm JHESS or via the recent algorithm JHMSH and its variants. 
 
 The main drawback of JHESS (or JHMSH) is that it may suffer from a fatal breakdown, causing a  brutal stop  of the computations and hence, the  $SR$-algorithm  does not run.  JHESS may also encounter  near-breakdowns, source of serious numerical instability.  
 
In this paper, we focus on these aspects. 
We first bring light on the necessary and sufficient condition for  the existence of the  $SR$-decomposition, which is intimately linked to $J$-Hessenberg reduction. Then we will 
  derive a strategy for curing  fatal breakdowns and also for treating  near breakdowns.  Hence,  
 the $J$-Hessenberg form  may be obtained.  
   Numerical experiments  are given,  demonstrating  the  efficiency of our strategies to  cure and treat breakdowns or near breakdowns.
\end{abstract}


\end{frontmatter}


\section{Introduction}

Let $A$ be a $2n\times 2n$ real matrix. The $SR$ factorization 
consists in writing $A$ as a product $SR$, where $S$ is
symplectic and $R=\left[
\begin{array}{ll}
R_{11} & R_{12}\\
 R_{21} & R_{22}
\end{array}
 \right]$
 is such that $R_{11},\;R_{12},\;R_{22}$ are upper triangular and
 $R_{21}$ is strictly upper triangular
  \cite{ Buns1, Dela2}. The factor $R$ is called $J$-triangular. This decomposition  plays an important role
 in structure-preserving methods for solving the
eigenproblem of a class of structured matrices.

More precisely, the $SR$ decomposition  can be interpreted  as the analog of the $QR$
decomposition \cite{Golu}, when instead of an Euclidean space, one considers a symplectic space : a 
linear space, equipped with a skew-symmetric inner product (see
for example \cite{Sal1} and the references therein). The 
 orthogonal group with respect to this indefinite inner product,  is  called the 
symplectic group and is unbounded  (contrasting with the Euclidean case).

In the literature,  the $SR$ decomposition is carried out, via the algorithm SRDECO, derived in \cite{Buns}.  SRDECO is based in the use  
of two kind of both symplectic and orthogonal transformations introduced in  \cite{Paig, Vloa}
and a third symplectic
but non-orthogonal  transformations,  proposed in
\cite{Buns}.  In fact, 
in \cite{Buns1}, it has been shown that $SR$ decomposition of 
a general matrix
can  not be performed  by employing 
only the above orthogonal and symplectic transformations.

We mention that the above transformations involved in SRDECO algorithm 
  are not  elementary rank-one modification of the identity (transvections), see \cite{Art, Golu}.

Recently in \cite{Sal2}, an algorithm, SRSH, based on symplectic transformations which are rank-one modification of the identity is derived, for computing the $SR$ decomposition.  These transformations   are called  symplectic Householder transformations. The new algorithm SRSH
     involves free parameters and advantages may be
taken from this fact.  An optimal version of SRSH, called SROSH is given in 
 \cite{Sal3}.
   Error analysis and computational aspects of this algorithm have been studied  
\cite{Sal4}. 

  In order to build a $SR$-algorithm (which is a $QR$-like algorithm) for computing the eigenvalues and eigenvectors of a matrix \cite{Wat}, a reduction of the matrix to an upper $J$-Hessenberg form is needed and is crucial.

In \cite{Buns}, the algorithm JHESS, for  reducing  a general matrix to an upper $J$-Hessenberg form is presented, using to this aim, an adaptation of SRDECO. 
In \cite{Sal6},  the algorithm JHSH, based on an adaptation of SRSH, is introduced, for 
 reducing a general matrix, to an upper $J$-Hessenberg form.  Variants of JHSH, named JHOSH, JHMSH and JHMSH2  are then derived, motivated by the numerical stability. 
 The algorithms JHESS as well as JHSH and its variants  have  $O(n^3)$ as complexity. 

The algorithm JHESS (and also SRDECO), as described in \cite{Buns}  may be subject of a fatal breakdown, causing a  brutal stop  of the computations. As consequence, the $J$-Hessenberg reduction can not be computed  and the  $SR$-algorithm  does not run. 
  Moreover, we demonstrate that the algorithm JHESS may breaks down while a condensed $J$ Hessenberg form  exists.  
These algorithms  also may suffer from severe form of near-breakdowns, source of serious numerical instability. 

In this paper, we restrict ourselves to the study  of such aspects, bringing significant insights on SRDECO and JHESS algorithms. 


 We will show derive a strategy for curing  fatal breakdowns   and treating near breakdowns. To this aim, 
     we first bring light on the $SR$-decomposition and SRDECO algorithm, in connection with 
     the theory developed by Elsner in \cite{Els}.
     Then,  a strategy for remedying to such breakdowns is proposed. 
      The same strategy is used for treating the near-breakdowns.   Numerical experiments  are given,  demonstrating  the  efficiency of our strategies to cure breakdowns or  to treat near breakdowns.


 The remainder of this paper is organized as follows. Section
 2, is devoted to the necessary preliminaries.   In the section 3,   the algorithms SRDECO, SRSH or SRMSH  are presented.
 Then, we  establish  a  connection between some coefficients of the current matrix produced by the SRDECO algorithm, when applied to a matrix $A$  and the necessary and sufficient condition for the existence of the $SR$ decomposition of $A$, as given in \cite{Els}. 
 In section 4, 
  we recall the algorithms JHESS and JHMSH. 
  We present then an example, for which 
 a fatal breakdown is  meet, in JHESS algorithm (also for JHMSH), 
  for reducing the matrix to an upper $J$-Hessenberg.  And hence, following the description of JHESS in \cite{Buns},  the algorithm stops. No hope to build then a SR-algorithm. 
We will show on the same example,  that such breakdown is curable. We then develop our strategies for  curing the fatal breakdowns. The same strategies are applied for near-breakdowns. The 
Section 5 is devoted to 
  numerical experiments and comparisons. We conclude in the section 6. 
  %
\section{Preliminaries}
 Let $J_{2n}$ (or simply $J$) be the $2n$-by-$2n$ real matrix
 \begin{equation}\label{Ji2n}
{ J}_{2n}=\left[\begin{array}{ll} %
{ 0}_n &{ I}_n\\ %
-{ I}_n & { 0}_n
\end{array}
\right],
\end{equation}
 where $0_n$ and $I_n$ stand respectively  for $n$-by-$n$ null and identity
 matrices.
The linear space  $\mathbb{R}^{2n}$  with the indefinite
skew-symmetric inner product
 \begin{equation}
 (x,y)_J = x^T J y
 \end{equation}
  is called symplectic. 
For 
$x,\;y \in\mathbb{R}^{2n},\mbox{ the orthogonality }x \perp' y$ stands for $(x,y)_J=0.$
   The symplectic adjoint $x^J$ of a vector  $x$,  is defined by
  \begin{equation}\label{deff}
  x^J = x^T J.
  \end{equation}
    The symplectic  adjoint of
    $M \in \mathbb{R}^{2n \times 2k}$ is defined  by
 \begin{equation}\label{def3}
 M^{J}=J_{2k}^T M^{T}J_{2n}.
 \end{equation}
 A matrix $S \in \mathbb{R}^{2n \times 2k}$ is called symplectic
 if
 \begin{equation}\label{def2}
 S^{J} S= I_{2k}.
 \end{equation}
 The symplectic  group (multiplicative group of square symplectic matrices) is denoted
 $\mathbb{S}.$
\noindent
A  transformation $T$  given by

\begin{equation}\label{sytr1}
 T=
I+ cvv^J \mbox{ where } c \in \mathbb{R},\; \; v \in
\mathbb{R}^{\nu}\;\; \mbox{ (with } \nu \mbox{ even),}
\end{equation}
 is called symplectic Householder transformation \cite{Sal2}. It satisfies
  \begin{equation}\label{jtra}
  T^J=I - c v
 v^J.
 \end{equation}
  The vector $v$ is called the direction of $T.$

For $x,\;y \in \mathbb{R}^{2n},\;$
  there exists  a symplectic Householder transformation $T$
   such that $Tx=y$
 if $x=y$ or $x^Jy\neq0.$ When $x^Jy\neq0,$
  $T$  is given by
 $$T=I-\frac{1}{x^J y}(y-x)(y-x)^J.$$
Moreover,
 each non null vector $x$ can be mapped onto any non null vector $y$
 by a product of at  most two symplectic Householder transformations \cite{Sal2}.
 Symplectic Householder transformations are rotations, i.e.
$det(T)=1$ and
 the  symplectic group $\mathbb{S}$ is generated by
 symplectic Householder transformations.
 In \cite{Paig, Vloa} two  orthogonal and symplectic transformations have  been introduced.
    The first, for which we refer as Van Loan's Householder transformation,   has the form
  \begin{equation}\label{reff1}
  H(k,w)=\left(
  \begin{array}{cc}
  \mbox{diag}\displaystyle{(I_{k-1},P)} & 0\\
  0 & \mbox{diag}\displaystyle{ (I_{k-1},P)}
  \end{array}
  \right),
  \end{equation}
    where
  $$P=I-2w w^T/w^Tw, \;\;w \in \mathbb{R}^{n-k+1}.$$
     The
  second, for which we refer as Van Loan's  Givens  transformation,  is
  \begin{equation}\label{reff2}
  J(k,\theta)=\left(
  \begin{array}{ll}
  C & S\\
  -S & C
  \end{array}
  \right),
\end{equation}
  where
  $$C=diag(I_{k-1},cos\theta,I_{n-k}) \mbox{ and }
S=diag(0_{k-1},sin\theta,0_{n-k}).$$
 $J(k,\theta)$ is a  Givens symplectic matrix, that is an
  "ordinary" 2$n$-by-$2n$ Givens rotation that
 rotates in planes $k$ and $k+n$ \cite{Wilk}.
%
The $SR$ factorization can not
     be performed for a general matrix by  using the sole
      $H(k,w)$ and $J(k,\theta)$
     transformations \cite{Buns}. A third type, introduced in \cite{Buns}, is given  by
  \begin{equation}
  G(k,\nu)=\left(\begin{array}{ll}
  D & F\\
  0 & D^{-1}
  \end{array}
  \right),
  \end{equation}
  where $k \in \{2,\ldots,n\},\,\nu \in \mathbb{R}$ and $D,\;
   F$ are the $n\times n$ matrices
   $$
   D=I_n + (\frac{1}{(1+{\nu}^2)^{1/4}}-1)(e_{k-1}e_{k-1}^T +
   e_{k}e_{k}^T),$$
   $$F= \frac{\nu}{(1+{\nu}^2)^{1/4}}(e_{k-1}e_{k}^T +
   e_{k}e_{k-1}^T).$$
   The matrix $G(k,\nu)$  is symplectic and non-orthogonal. 
   The  SRDECO algorithm  is then derived for computing    $SR$
   factorization for a general matrix, based on $ H, J$  and  $G$ transformations. 
A reduction of a general matrix to an upper $J$-Hessenberg, is obtained by using the same transformations involved in SRDECO, giving rise to JHESS algorithm. The breakdown in $SR$-decompoisition via SRDECO or in the reduction to an upper $J$-Hessenberg form via JHESS, when it occurs, is caused by the latest transformations $G.$ 
\section{SRDECO, SRSH algorithms}
The aim of this work, is to bring significant contributions on the understanding and the behaviour of the algorithms SRDECO, SRSH, JHESS and JHMSH.  Also, we propose  strategies for curing breakdowns. Similar strategies are applied also for remedying to near breakdowns.
\subsection{$SR$ decomposition : SRDECO, SRSH algorithms}
We consider the  SRDECO algorithm, as introduced in \cite{Buns}. Given $A \in \mathbb{R}^{2n\times 2n}$, the algorithm determines an $SR$ decomposition of $A$, using functions vlg, vlh and gal below. 
 The  function vlg uses  Van Loan's Givens transformation $J(k,c,s) $ as follows : for a given integer $1\leq k\leq n$ and a vector $a \in \mathbb{R}^{2n}$, it determines coefficients $c$ and $s$ such that the $n+k$th component of $J(k,c,s)a$ is zero. All components of $J(k,c,s)a$ remain unchanged, except eventually the $k$th and the $n+k$.
 \begin{alg}\label{vlg}
 	function[c,s]=vlg(k,a)\\
 	\indent $twon=length(a);\;n=twon/2;$\\
 	\indent $r=\sqrt{a(k)^2+a(n+k)^2};$\\
 	\indent if $r=0$ then $c=1;\;s=0;$\\
 	\indent else
 	\indent $\displaystyle{c=\frac{a(k)}{r};\;\;s=\frac{a{(n+k)}}{r}};$\\
 	\indent end
 \end{alg}
 The  function $vlh$ uses Van Loan's Householder transformation $H(k,w)$ as follows : for a given integer $k \leq n$ and a vector $a \in \mathbb{R}^{2n}$, a vector $w=(w_1,\ldots,w_{n-k+1})^T$ is determined such that the components $k+1,\ldots,n$ of $H(k,w)a$ are zeros. All components $1,\ldots,k-1$ and $n+1,\ldots,n+k-1$ remain unchanged.
 \begin{alg}\label{vlh}
 	function[$\beta$,w]=vlh(k,a)\\
 	\indent $twon=length(a);\;n=twon/2;$\\
 	\indent \% $w=(w_1,\ldots,w_{n-k+1})^T;$\\
 	\indent  $r1=\sum_{i=2}^{n-k+1} a(i+k-1)^2;$\\
 	\indent  $r=\sqrt{a(k)^2+r1};$\\
 	\indent $w_1= a(k)+sign(a(k))r;$\\
 	\indent $w_i=a{(i+k-1)}$ for $i=2,\ldots,n-k+1;$\\
 	\indent  $r={w_1^2+r1};\;\;
 	\displaystyle{\beta=\frac{2}{r};}$\\
 	\indent \%$ P=I-\beta ww^T;\;\;(H(k,w)a)_i=0$ for $i=k+1,\ldots,n.$\\
 	\indent end
 \end{alg}
 The  function $gal$ uses the  transformation $G(k,\nu)$ as follows : for a given integer $k \leq n$ and a vector $a \in \mathbb{R}^{2n}$, satisfying the condition $a_{n+k}=0$ only if $a_{k+1}=0$, it determines $\nu$ such that the $k+1$th of $G(k,\nu)a$ is zero. 
 \begin{alg}\label{gal}
 	function$[\nu]=gal(k,a)$\\
 	\indent $twon=length(a);\;n=twon/2;$\\
 	\indent if $a_k=0$\\
 	\indent $\nu=0;$\\
 	\indent else\\
 	\indent $\nu=-\displaystyle{\frac{a_{k+1}}{a_{n+k}}}$;\\
 	\indent end\\
 	\indent end
 \end{alg}


The  algorithm SRDECO is as follows : the matrix $A$ is overwritten by the $J$-upper triangular matrix. If $A$ has no $SR$ decomposition, the algorithm stops.
   \begin{alg}
   	function [S,A]=SRDECO(A)\\
   \textbf{1.}	For $j=1,\ldots,n$\\
  \textbf{2.} 	For $k=n,\ldots,j$\\
   \textbf{3.}	    Zero the entry $(n+k,j)$  of $A$ by running the function $[c,s]=vlg(k,A(:,j))$ 
  and computing  
   	    $J_{k,j}=J(k,c,s).$\\
   	  \textbf{4.}    Update $A=J_{k,j}A$ and $S=S J_{k,j}^T$\\
 \textbf{5.}  	    End for.\\
  \textbf{6.} 	    Zero the entries $(j+1,j),\dots,(n,j)$ of $A$ by running the function $\displaystyle{[\beta,w]=vlh(j,A(:,j))}$ 
    and computing $H_j=H(j,w).$ \\
      \textbf{7.}  Update $A=H_jA$ and $S=S H_j^T$.\\
   \textbf{8.}	    If  $j\leq n-1$\\
   \textbf{9.}	     For $k=n,\ldots,j+1$\\
   	\textbf{10.}     Zero the entry $(n+k,n+j)$ of $A$ by running the function $\displaystyle{[c,s]=vlg(k,A(:,n+j))}$ 
  and computing  $J_{k,n+j}=J(k,c,s).$ \\
   	  	\textbf{11.}     Update $A=J_{k,n+j}A$ and $S=S J_{k,n+j}^T.$\\
   	 \textbf{12.}    End for.\\
   	  \textbf{13.}   Zero the entries $(j+2,n+j),\dots, (n,n+j)$ of $A$ by running the  
   	function $\displaystyle{[\beta,w]=vlh(j+1,A(:,j))}$ and computing $H_{n+j}=H(j+1,w)$. \\
   	  \textbf{14.} Update $A=H_{n+j}A$ and $S=S H_{n+j}^T.$\\
   	   \textbf{15.}  If the entry $(j+1,n+j)$ of $A$ is nonzero and the entry $(n+j,n+j)$ is 
   	    zero then stop the algorithm,\\
   	   \textbf{16.}  else\\
   	    \textbf{17.}  Zero the entry $(j+1,n+j)$ of $A$ by running the function $\displaystyle{[\nu]=gal(j+1,A(:,n+j))}$ and computing 
   	     $G_{j+1}=G(j+1,\nu).$ \\
   	       \textbf{18.} Update $A=G_{j+1}A$ and $S=S G_{j+1}^{-1}.$ \% $G_{j+1}^{-1}=G_{j+1}^{J}.$\\
   	    \textbf{19.} End if\\
   	    \textbf{20.} End if\\
   	    \textbf{21.} End for.
   	       \end{alg}
   	    The $SR$ decomposition can be also performed using only the symplectic Householder transformations $T$ of (\ref{sytr1}), which are rank-one modifications of the identity, giving rise to the algorithm SRSH. More on this can be found in \cite{Sal2, Sal3, Sal4}. A modified version of SRSH, numerically more stable,  is SRMSH. It turn out that SRMSH shares the same steps (1-16) of SRDECO, but not the remaining ones. In fact, the function $gal$ and the symplectic matrices $G_{j+1}$ in SRDECO are replaced by the function 
   	    $sh2$ and $T_j$ having the form of (\ref{sytr1}). The function $sh2$ is as follows
   	    \begin{alg}\label{sh2}
   	    	function [c, v] = sh2(a)\\
   	    	\indent \%compute $c$ and $v$ such that $T_2 e_1 = e_1,$ and $T_2 a = \mu e_1 + \nu e_{n+1},$ \\
   	    	\indent \%$\mu$ is a free parameter, and  $T_2=(eye(twon)+c * v * v' * J);$\\
   	    	\indent $twon=length(a);\;n=twon/2;$\\
   	    	\indent $J=[zeros(n),eye(n);-eye(n),zeros(n)];$\\
   	    	\indent If $n==1$\\
   	    	\indent\indent    $v=zeros(twon,1);\;c=0;$   $\%T=eye(twon);$\\
   	    	\indent else\\
   	    	\indent Choose $\mu;$\\
   	    	\indent       $ \nu =a(n+1);$\\
   	    	\indent        If $ \nu ==0$\\
   	    	\indent           Display('division by zero')\\
   	    	\indent           Return\\
   	    	\indent       else\\
   	    	\indent           $v=\mu e_1 + \nu e_{n+1} - a,$ 
   	    	${\displaystyle c=\frac{1}{a(n+1)(a(1)-\mu)};}$\\
   	    	%
   	    	\indent  End\\
   	    	\indent End
   	    \end{alg}
   	  We obtain the algorithm
   	     \begin{alg}
   	     	function [S,A]=SRMSH(A)\\
   	    \textbf{1.} Run steps	\textbf{1.}-\textbf{16.} of SRDECO\\     
   	     \textbf{17.}  Set $c0=[j:n,n+j:2n],\;c1=[j+1:n,n+j:2n]$, \\
   	      \textbf{18.} Zero the entry $(j+1,n+j)$ of $A$ by running the function $\displaystyle{[c,v]=sh2(A(co,n+j))}$  \\
   	       \textbf{19.} compute  
   	     $T_{j}=I+c vv^TJ,$ update $A(:,c1)=T_{j}A(:,c1)$ and 
   	      $S(:,c0)=S T_{j}^{J}.$\\
   	     \textbf{20.} End if\\
   	     \textbf{21.} End if\\
   	     \textbf{22.} End for.
   	   \end{alg} 
   	 \begin{rmk}
   	 	The function $sh2$ in the body of the algorithm $SRMSH$ may by replaced by the function $osh2$ (see \cite{Sal3, Sal4}) which presents the best conditioning among all possible choices. 
   	 	\end{rmk}   
\subsection{Discussion : existence of SR decomposition, link with SRDECO and SRMSH}  
In this subsection, we bring light on the connection between the existence of  $SR$ decomposition and the algorithm SRDECO or equivalently SRMSH. We recall first the following result, given in \cite{Els} :
\begin{thm}\label{Elsn}
	Let $A \in \mathbb{R}^{2n\times 2n}$ be nonsingular and $P$ the permutation matrix $P=[e_{1},e_{n+1},e_{2},e_{n+2},\ldots,e_{n},e_{2n}]$, where $e_i$ denotes the $i$th canonical vector of $\mathbb{R}^{2n}.$ There exists  $S\in \mathbb{R}^{2n\times 2n}$ symplectic and  $R\in \mathbb{R}^{2n\times 2n}$ upper $J$-triangular, such that $A=SR$ if and only if all even leading minors of $P^T A^T J A P$ are nonzero. 
\end{thm}
In \cite{Buns}, a comment on SRDECO states : "if at any stage $j\in \{1,\ldots,n-1\}$ the algorithm ends because of the stopping condition, then the $2j$th leading principal minor of $P^T A^T J A P$ is zero, and $A$ has no $SR$ decomposition (see Theorem \ref{Elsn})."
However, a proof of how is  connected   the stopping condition of the algorithm SRDECO with the condition of  Theorem \ref{Elsn} is not given. Remark also that for SRDECO algorithm, the condition $A$ nonsingular is not required, while it is  for Theorem \ref{Elsn}. Here we give a proof on how this connection is made. 
Notice first that if $A=SR$, where $S$ is any symplectic matrix and $R$ is any matrix, then $A^TJ A=R^T S^T JSR=R^TJR$. Hence a minor of $A^TJA$ is equal to its corresponding one of $R^TJR$. The same equality between minors is valid also for $P^T A^T JAP$ and $P^T R^T JRP.$
The following Theorem establishes an explicit relation between the leading $2j$-by-$2j$ minors of $P^T A^T JAP$ and the computed coefficients which determine  the stopping condition of SRDECO. For a given matrix $M$, let us denote by $M_{[j,j]}$ the submatrix obtained from $M$ by deleting all rows and columns except rows and columns $1,\ldots,j.$ We have
\begin{thm}\label{con}
 Let  $A \in \mathbb{R}^{2n\times 2n}$ be a  matrix (not necessarily nonsingular), and let $R$ be the matrix that one  obtains  at stage $1 \leq j \leq n-1$ of the algorithm SRDECO, by executing instructions {\bfseries{1.}} to {\bfseries{14.}} (corresponding to the current updated  matrix $A$ in the process, at stage $j$ and until instruction {\bfseries{14.}}).  Then the leading $2j$-by-$2j$ minor of $P^T A^T JAP$ satisfies
 \begin{equation}
 det((P^T A^T J A P)_{[2j,2j]}) =\left[r_{1,1}r_{n+1,n+1}\ldots r_{i,i}r_{n+i,n+i}\ldots r_{j,j}r_{n+j,n+j}\right]^2 .
 \end{equation}
 \end{thm}
 \begin{pr}
 Partitioning 
 	$R=\left(
 	\begin{array}{ll}
 	R_{11} & 	R_{12} \\
 		R_{21} & 	R_{22} \\
 	\end{array}
 	\right),$ then $R_{11},\; 	R_{12},\; 	R_{21},\; 	R_{22}$ have the form  \\
 	$$R_{11}=
 	 \left[
 	 \begin{array}{lllllll}
 	 r_{1,1}&r_{1,2}& \ldots & r_{1,j} & r_{1,j+1} & \ldots & r_{1,n} \\
 	 0&r_{2,2}& \ldots & r_{2,j} & r_{2,j+1} & \ldots & r_{2,n} \\
 	 \vdots &\ddots& \ddots & \vdots & \vdots& \ldots & \vdots \\
 	 \vdots && \ddots & r_{j,j} & r_{j,j+1} & \ldots & r_{j,n} \\
 	 0 &\ldots& \ldots & 0 & r_{j+1,j+1} & \ldots & r_{j+1,n}     \\
 	 \vdots &&  & \vdots & \vdots&  & \vdots   \\
 	 0 && \ldots & 0 & r_{n,j+1} & \ldots & r_{n,n} 
 	 \end{array}
 	 \right],
 	$$
 	$$R_{12}=
 	\left[
 	\begin{array}{llllllll}
 	r_{1,n+1}&r_{1,n+2}& \ldots & r_{1,n+j-1} & r_{1,n+j} & r_{1,n+j+1} &\ldots & r_{1,2n} \\
 	0&r_{2,n+2}& \ldots & r_{2,n+j-1} & r_{2,n+j} & r_{2,n+j+1} & \ldots & r_{2,2n} \\
 	\vdots &\ddots& \ddots & \vdots & \vdots& \vdots&  & \vdots \\
 	\vdots && \ddots & r_{j-1,n+j-1} & r_{j-1,n+j} & r_{j-1,n+j+1} & \ldots & r_{j-1,2n} \\
 	0 &\ldots& \ldots & 0 & r_{j,n+j} & r_{j,n+j+1} & \ldots & r_{j,2n}     \\
 	 	0 &\ldots& \ldots & 0 & r_{j+1,n+j} & r_{j+1,n+j+1} & \ldots & r_{j+1,2n}     \\
 	 	 	 	0 &\ldots& \ldots & 0 & 0 & r_{j+2,n+j+1} & \ldots & r_{j+2,2n}     \\
 	\vdots &&  & \vdots & \vdots&\vdots&  & \vdots   \\
 	0 && \ldots & 0 & 0 & r_{n,n+j+1} & \ldots & r_{n,2n} 
 	\end{array}
 	\right],
 	$$
  $$R_{21}=\left[
  \begin{array}{lllllll}
  0&   r_{n+1,2}& \ldots & r_{n+1,j} & r_{n+1,j+1} & \ldots & r_{n+1,n}\\
  \vdots &\ddots& \ddots & \vdots & \vdots&   & \vdots \\
  0&\ldots&  0& r_{n+j-1,j} & r_{n+j-1,j+1} & \ldots & r_{n+j-1,n} \\
  0 &\ldots&  0& 0 & r_{n+j,j+1} & \ldots & r_{n+j,n}\\
  \vdots &&  \vdots & \vdots & \vdots&  & \vdots\\
  0 &\ldots& 0 & 0 & r_{2n,j+1} & \ldots & r_{2n,n} 
  \end{array}
  \right],$$ and
  $$R_{22}=
  \left[
  \begin{array}{lllllll}
  r_{n+1,n+1}&r_{n+1,n+2}& \ldots  & r_{n+1,n+j} & r_{n+1,n+j+1} &\ldots & r_{n+1,2n} \\
  0&r_{n+2,n+2}& \ldots  & r_{n+2,n+j} & r_{n+2,n+j+1} & \ldots & r_{n+2,2n} \\
  \vdots &\ddots& \ddots  & \vdots& \vdots&  & \vdots \\
  \vdots&\ldots& \ddots  & r_{n+j,n+j} & r_{n+j,n+j+1} & \ldots & r_{n+j,2n}     \\
  0 &\ldots& \ldots & 0 & r_{n+j+1,n+j+1} & \ldots & r_{n+j+1,2n}     \\
  \vdots &&  & \vdots &\vdots&  & \vdots   \\
  0 &\ldots& \ldots & 0 & r_{2n,n+j+1} & \ldots & r_{2n,2n} 
  \end{array}
  \right].
  $$
  The stopping condition at this stage $j$ is "$r_{n+j,n+j} = 0$ and $r_{j+1,n+j} \neq 0$".
   We will establish connection between  the coefficients $r_{i,i},\;r_{n+i,n+i}$ of the current matrix $R$, with $1\leq i \leq j$ and the leading $2j$-by-$2j$ minor of $P^T A^T J A P.$  Setting 
   $\hat{J}=P^TJ P$ and $\hat{R}=P^TR P$, we get 
   $P^T A^T J A P=P^T R^T J R P=(P^T R^T P)(P^TJ P )(P^TR P)=\hat{R}^T \hat{J}\hat{R}.$ Recall that $\hat{J}=diag(J_2,\ldots,J_2).$  Partitioning 
   $\hat{R}=\left(
   \begin{array}{ll}
   \hat{R}_{11} & 	\hat{R}_{12} \\
   \hat{R}_{21} & 	\hat{R}_{22} \\
   \end{array}
   \right),$ with the block $\hat{R}_{11}$ is $2j$-by-$2j$. Then we obtain for $\hat{R}_{11},\; 	\hat{R}_{12},\; 	\hat{R}_{21},\; 	\hat{R}_{22}$ :   
   \begin{equation}\label{upper}
   \hat{R}_{11}=
   \left[
   \begin{array}{llllllll}
   r_{1,1}&r_{1,n+1}& \ldots & r_{1,i} & r_{1,n+i} & \ldots & r_{1,j} & r_{1,n+j} \\
     0&r_{n+1,n+1}& \ldots & r_{n+1,i} & r_{n+1,n+i} & \ldots & r_{n+1,j} & r_{n+1,n+j} \\
     \vdots & 0&\ddots &\vdots &\vdots &\vdots &\vdots &\vdots \\
          \vdots & \vdots  & \ddots & r_{i,i} & r_{i,n+i} & \ldots & r_{i,j} & r_{i,n+j} \\
                  \vdots &  \vdots & & 0& r_{n+i,n+i} & \ldots & r_{n+i,j} & r_{n+i,n+j} \\
                 \vdots &  \vdots &  & 0& 0 & \ddots & \vdots & \vdots \\ 
                     \vdots & \vdots  & & 0 & 0 & \ddots & r_{j,j} & r_{j,n+j} \\
                       0&  0 & \ldots & 0 & 0 & \ldots & 0 & r_{n+j,n+j} 
   \end{array}
   \right],
   \end{equation}
    which is a upper $2j$-by-$2j$ triangular matrix. The $2(n-j)$-by-$2j$  block  $\hat{R}_{21}$ turn out to have all entries zeros except the entry in position $(1,j).$ More precisely,
     $\hat{R}_{21} = \left(
     \begin{array}{ll}
     0 & r_{j+1,n+j}\\
     0  & 0
     \end{array}
     \right).
     $
     Setting $\hat{J}_{2k}= diag( J_2,\ldots,J_2)  \in  \mathbb{R}^{2k \times 2k} $ for an integer $k$,
      and due to the special structures of $\hat{J}$, we get \\
      \begin{tabular}{lll}
      	 $\hat{R}^T \hat{J} \hat{R}$ & = & $\left(
      	 \begin{array}{ll}
      	 	\hat{R}_{11} & 	\hat{R}_{12} \\
      	 	\hat{R}_{21} & 	\hat{R}_{22} \\
      	 \end{array}
      	 \right)^T 
      \hat{J}	 
      \left(
      \begin{array}{ll}
      \hat{R}_{11} & 	\hat{R}_{12} \\
      \hat{R}_{21} & 	\hat{R}_{22} \\
      \end{array}
      \right)	 
      	 $ \\
      	   & = & $\left(
      	  \begin{array}{ll}
      	  \hat{R}_{11}^T & 	\hat{R}_{21}^T \\
      	  \hat{R}_{12}^T & 	\hat{R}_{22}^T \\
      	  \end{array}
      	  \right)
      	  \left(
      	  \begin{array}{ll}
      	 \hat{J}_{2j}	  \hat{R}_{11} & \hat{J}_{2j}	 	\hat{R}_{12} \\
      	\hat{J}_{2(n-j)}	   \hat{R}_{21} & \hat{J}_{2(n-j)}	 	\hat{R}_{22} \\
      	  \end{array}
      	  \right).$	 
      	 	\end{tabular} 
      	 	
      	 	 Let $(\hat{R}^T \hat{J} \hat{R})_{[2j,2j]}$ denote the leading $2j$-by-$2j$ block of $\hat{R}^T \hat{J} \hat{R},$ we obtain 
      	 	 $$(\hat{R}^T \hat{J} \hat{R})_{[2j,2j]}=  \hat{R}_{11}^T   \hat{J}_{2j}	  \hat{R}_{11} + 
      	 	 	\hat{R}_{21}^T    	\hat{J}_{2(n-j)}	   \hat{R}_{21} .$$
      	 	 	Denoting $\tilde{e}_j=(0,\ldots,0,1)^T$ the $j$th canonical vector of $\mathbb{R}^j$ and 
      	 	 	$e_1=(1,0,\ldots,0)^T,$  $e_2=(0,1,0,\ldots,0)^T$ respectively the first and the second canonical vectors of $\mathbb{R}^{2(n-j)},$ then $\hat{R}_{21}$ may be expressed as 
      	 	 	 $\hat{R}_{21}= r_{j+1,n+j} e_1 \tilde{e}_{2j}^T.$ Hence, 
      	 	 	  $\hat{R}_{21}^T    	\hat{J}_{2(n-j)}	   \hat{R}_{21} = r_{j+1,n+j}^2  \tilde{e}_{2j}e_1^T 
      	 	 	   \hat{J}_{2(n-j)}e_1 \tilde{e}_{2j}^T.$ As $  \hat{J}_{2(n-j)}e_1=-e_2$,  we obtain directly 
      	 	 	   $\hat{R}_{21}^T    	\hat{J}_{2(n-j)}	   \hat{R}_{21} = -r_{j+1,n+j}^2  \tilde{e}_{2j}e_1^T 
      	 	 	   e_2 \tilde{e}_{2j}^T=0.$ Thus
      	 	 	    $$(\hat{R}^T \hat{J} \hat{R})_{[2j,2j]}=  \hat{R}_{11}^T   \hat{J}_{2j}	  \hat{R}_{11} .$$
      	 	 	    It follows 
      	 	 	    $$det((\hat{R}^T \hat{J} \hat{R})_{[2j,2j]})=  det(\hat{R}_{11}^T) det(  \hat{J}_{2j}	)  det(\hat{R}_{11}) .$$
      	 	 	    As $det(  \hat{J}_{2j}	 )=1$, and    $det(\hat{R}_{11}^T)=det(\hat{R}_{11}) ),$  we get 
      	 	 	      $$det((\hat{R}^T \hat{J} \hat{R})_{[2j,2j]})=  (det(\hat{R}_{11}))^2 .$$
      	 	 	    
      	 	 	    Since $P^T A^T J A P =(\hat{R}^T \hat{J} \hat{R}),$ it follows 
      	 	 	    that $$(P^T A^T J A P)_{[2j,2j]} =(\hat{R}^T \hat{J} \hat{R})_{[2j,2j]} ,$$
      	 	 	    which implies for the $2j$-by-$2j$ leading minor of $P^T A^T J A P$
      	 	 	    $$det((P^T A^T J A P)_{[2j,2j]}) =det((\hat{R}^T \hat{J} \hat{R})_{[2j,2j]})=  (det(\hat{R}_{11}))^2.$$
      	 	 	    The matrix $\hat{R}_{11}$ is $2j$-by-$2j$ upper triangular matrix, and from relation (\ref{upper}),  we have 
      	 	 	    $$det((P^T A^T J A P)_{[2j,2j]}) =\left[r_{1,1}r_{n+1,n+1}\ldots r_{i,i}r_{n+i,n+i}\ldots r_{j,j}r_{n+j,n+j}\right]^2 .$$
      	 	 	  \end{pr}
      	 	 	  \begin{cor}
      	 	 	  	Let  $A \in \mathbb{R}^{2n\times 2n}$ be a nonsingular matrix, 
      	 	 	and let $R$ be the matrix that one  obtains  at stage $1 \leq j \leq n-1$ of the algorithm SRDECO, by executing instructions {\bfseries{1.}} to {\bfseries{14.}} (corresponding to the current updated  matrix $A$ in the process, at stage $j$ and until instruction {\bfseries{14.}}).  	
      	 	 	  	Then 
      	 	 	  	 $A$ admits an $SR$ decomposition if and only if $r_{n+j,n+j} \neq 0,\;\forall j \in\{1,\ldots,n\}.$
      	 	 	  \end{cor}
      	 	 	  \begin{pr}
      	 	 	  	Since $A$ is nonsingular and using Threorem \ref{Elsn} and Theorem \ref{con}, we have :  $A$ admits a $SR$ decomposition if and only if $r_{1,1}r_{n+1,n+1}\ldots
      	 	 	  	 r_{j,j}r_{n+j,n+j} \neq 0, \forall j \in\{1,\ldots,n\}.$ At the stage $j$,  we have $A=SR$ for some symplectic matrix $S.$  Due to the structure of $R$, we deduce that the coefficients $r_{1,1},r_{2,2},\ldots,
      	 	 	  	 r_{j,j}$ are automatically all nonzero (otherwise $R$ would be singular and so would be $A$). The result is then straightforward.
      	 	 	  	\end{pr}
     	If the condition $A$ nonsingular is not required,  one may ask in this case whether the $SR$-decomposition exits even when a $2j$-by-$2j$ leading minor $det((P^T A^T J A P)_{[2j,2j]})$  is equal zero for some $j.$  We precise this in the following  result  
     	  \begin{thm}
     	  	Let  $A \in \mathbb{R}^{2n\times 2n}$ be any matrix, 
     	  	and let $R$ be the matrix that one  obtains  at stage $1 \leq j \leq n-1$ of the algorithm SRDECO (or SRMSH), by executing instructions {\bfseries{1.}} to {\bfseries{14.}} (corresponding to the current updated  matrix $A$ in the process, at stage $j$ and until instruction {\bfseries{14.}})
     	  Then 
     	  	$A$ admits an $SR$ decomposition if and only if 
     	  	$(r_{n+j,n+j}  \neq 0  \mbox{ or } r_{j+1,n+j}  =0),\;\forall j \in\{1,\ldots,n \}.$
     	  \end{thm}	 
     	  \begin{pr}
     	  	The condition is sufficient, since if it is satisfied, the stopping condition in SRDECO (or SRMSH) is never meet and a $SR$ decomposition is furnished at the end of the process. We show now that the condition is necessary, i.e. we show that if there exists an index $j$ such that 
     	  	 $r_{n+j,n+j}= 0$ and $r_{j+1,n+j}\neq 0$, then  $SR$ decomposition does not exist. In the fact, suppose that there exists an integer $1 \leq j \leq n-1$ such that 
     	  	  $r_{n+j,n+j}= 0$ and $r_{j+1,n+j}\neq 0$ and  let us seek for a symplectic matrix $S_j$ such that the product $S_j a=a$ for any vector $a$ possessing  the same structure of  any  column $1,\ldots,j$ and $n+1,\ldots,n+j-1$ of $R$ and transforms the $n+j$ th column $R(:,n+j)= \sum_{i=1} ^{j+1}r_{i,n+j} e_i + \sum_{i=1} ^{j-1}r_{n+i,n+j} e_{n+i}$ into the desired form 
     	  	 \begin{equation}\label{dfr1}
     	  	 S_j R(:,n+j)= \sum_{i=1} ^{j}r'_{i,n+j} e_i + \sum_{i=1} ^{j}r'_{n+i,n+j} e_{n+i}.
     	  	 \end{equation}  The matrix $S_j$ has necessarily the form
     	  	 $$S_j=[e_1,\ldots,e_j, s_{j+1},\ldots,s_{n},e_{n+1},\ldots,e_{n+j-1}, s_{n+j},\ldots,s_{2n}],$$
     	  	  where $e_k$ stands for the $k$th canonical vector of $\mathbb{R}^{2n}.$ Hence we get 
     	  	   	 \begin{equation}\label{dfr2}
     	  	  S_j R(:,n+j)= \sum_{i=1} ^{j}r_{i,n+j} e_i + 
     	  	  r_{j+1,n+j} s_{j+1}+
     	  	  \sum_{i=1} ^{j-1}r_{n+i,n+j} e_{n+i} .
     	  	  \end{equation}
     	  	  In one hand, from relation (\ref{dfr1}), we get $e_j^T S_j R(:,n+j)= r'_{n+j,n+j}.$ In the other hand, from relation (\ref{dfr2}), and the fact that $S_j$ is symplectic, we get 
     	  	  $e_j^T S_j R(:,n+j)= 0.$  Thus, we deduce  $r'_{n+j,n+j}=0.$ Therefore, the relations 
     	  	  (\ref{dfr1} - \ref{dfr2}), imply $r_{j+1,n+j}s_{j+1}$ belongs to the space spanned by $\{e_1,\ldots,e_j, e_{n+1},\ldots,e_{n+j-1}\}$.  Since the vectors of 
     	  	  $e_1,\ldots,e_j, s_{j+1},e_{n+1},\ldots,e_{n+j-1}$ are linearly independents,  we deduce $r_{j+1,n+j}=0,$   which is absurd. The matrix $S_j$ does not exist and hence $SR$ decomposition does not exist.
     	  	\end{pr}	    
  \begin{rmk}
  	Remark that $S_j$ corresponds to the symplectic matrix $G_{j+1}$ for SRDECO and to the symplectic matrix $T_j$ for SRMSH.
  	\end{rmk}
\section{Curing breakdowns or treating near-breakdowns in JHESS, JHMSH algorithms}
\subsection{ Breakdowns or near-breakdowns in JHESS, JHMSH algorithms}
The algorithm SRDECO  may be adapted for reducing a matrix to the condensed 
 upper $J$-Hessenberg form, see \cite{Buns}. This leads to the algorithm JHESS. In a similar way, the algorithm SRSH or its variant SRMSH may be adapted for handling the reduction of a matrix to $J$-Hessenberg form, see \cite{Sal6}. 
We recall that a matrix 
$H=
\left[
\begin{array}{ll}
H_{11} &H_{12} \\
H_{21} & H_{22}
\end{array}
\right]
\in  \mathbb{R}^{2n\times 2n},$ is upper $J$-Hessenberg when  $H_{11},\,H_{21},\,H_{22}$ are upper triangular and $H_{12}$ is upper Hessenberg. $H$ is called unreduced when $H_{21}$ is nonsingular and the Hessenberg $H_{12}$ is unreduced, i.e. the entries  of the subdiagonal are all nonzero.

The algorithm JHESS is formulated in \cite{Buns} as follows : "given a matrix $A\in \mathbb{R}^{2n\times 2n}$ and $S=I_{2n}$, the following algorithm reduces, if it is possible, $A$ to upper $J$-Hessenberg form $H=\Pi^{-1}A\Pi$, with a symplectic matrix $\Pi$ whose first column is a multiple of $e_1$. A is overwritten by the $J$-Hessenberg matrix $H$ and $S$ is overwritten by the transforming matrix $\Pi.$ If this reduction of $A$ does not exist, the algorithm stops".
 \begin{alg}
 	function [S,A]=JHESS(A)\\
 	\textbf{1.}	For $j=1,\ldots,n-1$\\
 	\textbf{2.} 	For $k=n,\ldots,j+1$\\
 	\textbf{3.}	    Zero the entry $(n+k,j)$  of $A$ by running the function $[c,s]=vlg(k,A(:,j))$ 
 and computing  
 	$J_{k,j}=J(k,c,s).$\\
 	\textbf{4.} 	Update $A=J_{k,j}AJ_{k,j}^T$ and $S=S J_{k,j}^T$\\
 	\textbf{5.}  	    End for.\\
 	\textbf{6.} 	    Zero the entries $(j+2,j),\dots,(n,j)$ of $A$ by running the function $\displaystyle{[\beta,w]=vlh(j+1,A(:,j))}$ 
   and computing $H_j=H(j+1,w).$ \\
   	\textbf{7.} Update $A=H_jAH_j^T$ and $S=S H_j^T$.\\
 		\textbf{8.}  If the entry $(j+1,j)$ of $A$ is nonzero and the entry $(n+j,j)$ is zero
 	 then stop the algorithm\\
 			\textbf{9.}  else\\
 			\textbf{10.}  Zero the entry $(j+1,j)$ of $A$ by running the function $\displaystyle{[\nu]=gal(A(:,j))}$ \\
 			\textbf{11.}  Compute 
 			$G_{j+1}=G(j+1,\nu).$ \\
 			\textbf{12.}  Update $A=G_{j+1}AG_{j+1}^{-1}$ and $S=S G_{j+1}^{-1}.$\\
 			\textbf{13.} End if\\
 			\textbf{14.}	     For $k=n,\ldots,j+1$\\
 			\textbf{15.}     Zero the entry $(n+k,n+j)$ of $A$ by running the function $\displaystyle{[c,s]=vlg(k,A(:,n+j))}$  and 
 		compute  $J_{k,n+j}=J(k,c,s).$ \\
 			\textbf{16.} 	Update $A=J_{k,n+j}AJ_{k,n+j}^T$ and $S=S J_{k,n+j}^T.$\\
 			\textbf{17.}    End for.\\
 	\textbf{18.}	    If  $j\leq n-2$\\
 	\textbf{19.}   Zero the entries $(j+2,n+j),\dots, (n,n+j)$ of $A$ by running the function
   $\displaystyle{[\beta,w]=vlh(j+1,A(:,n+j))}$ and compute $H_{n+j}=H(j+1,w)$. \\
 	\textbf{20.}   Update $A=H_{n+j}AH_{n+j}^T$ and $S=S H_{n+j}^T.$\\
    \textbf{21.}  End if\\
 	\textbf{22.}  End for.
 \end{alg}
One of the main drawback of JHESS is that a fatal breakdown can be encountered. To illustrate our purpose, we consider the following example. Let $A_6$ be the 6-by-6 matrix 
\begin{equation}\label{exam}
A_6=\left(
\begin{array}{llllll}
1  & 0   &0 &  1 &  2 &  0\\
2  & 1 &  0 &  2 &  1 &  0\\
0  & 2 &  1 &  0  & 2 &  1\\
0   &2  & 0  & 1 &  0 &  0\\
0  & 1  & 2  & 3 &  1 &  0\\
0 &  0  & 1  & 0  & 3  & 1
\end{array}
\right).
\end{equation}
 The algorithm JHESS, applied to $A_6$, meets a fatal breakdown at the first step : the entry $A_6(2,1)\neq 0$ and the entry $A_6(4,1)=0,$ the algorithm stops. In fact, it is impossible to find a symplectic matrix $S_1$, with the first column proportional to $e_1$ such that $SA_6e_1= \alpha e_1 + \beta e_{4}$, as showed in the above subsection.  Thus, $A_6$  can not be reduced to an upper $J$-Hessenberg form, via symplectic similarity transformations, for which the first column is proportional to $e_1.$  The SR-algorithm as described in \cite{Buns} , uses first   JHESS algorithm for reducing  a matrix to the $J$-Hessenberg form. As consequence, if applied to $A_6$, the SR-algorithm stops also at the first step.  Let us remark also that  the basic SR-algorithm (which can be roughly described as consisting in repeating the factorisation $A=SR$, and the product $A=RS$) works  and converges, when applied to the example $A_6.$  The algorithm JHESS may also suffers from another serious problem : the near-breakdown. The latter occurs when the condition number of the symplectic and non-orthogonal  matrix $G_{j+1}$ at the step 	\textbf{11.} of the algorithm JHESS becomes very large.  This causes a dramatic growth of the rounding errors. 
 
 The following strategy is proposed in \cite{Buns} for remedying to such problems : if in the $j$th iteration, the condition number of the matrix $G_{j+1}$ is larger than a certain tolerance, the iteration is stopped. In the implicit form (which is the useful one) of the algorithm SR, an exceptional similarity transformation is computed, with the symplectic (but non-orthogonal) matrix $S_j= I -ww^TJ$, where $w$ is a random vector with $\norm{w}_2=1.$ The algorithm JHESS is then applied to the new similar matrix $S_j^{-1}AS_j.$ If the number of encountered near-breakdowns/breakdowns exceeds a given bound, the whole process is definitively stopped. This strategy presents certain serious drawbacks : 1)  The condition number of $S_j^{-1}AS_j$ will be worse than the condition number of $A.$ This du to the fact that $S_j$ can never be orthogonal. Hence, numerical instability is expected. 2) The cost of forming the product $S_j^{-1}AS_j$  is $O(n^2)$ where $2n$ is the dimension of $A.$  3) The product $S_j^{-1}AS_j$ fills-up the matrix and destroys the previous partially created $J$-Hessenberg form of $A.$ Hence an additional cost of $O(n^3)$ is needed to restore the $J$-Hessenberg form.  To summarize, each application of this strategy creates a current matrix with worse condition number than the previous, and needs an expensive cost of  $O(n^3)$ flops. 
 
 In the sequel, we propose two  alternatives, for which  either all or some of the above drawbacks are avoided. The first consists in a careful choice of the random vector $w$ so that the product $S_j^{-1}AS_j$ does not fill-up the matrix and preserves all the created zeros during the previous steps $1,\ldots,j-1.$ This diminish considerably the cost. However, the condition number of $S_j^{-1}AS_j$ may become worse than this of $A.$ The second alternative is more attractive since it allows us to avoid all of the above drawbacks. It consists in computing  a similarity transformation $S_j^{-1}AS_j$  for which : 1) the proposed matrix $S_j$ is not only symplectic but also orthogonal. Thus, the condition number of $S_j^{-1}AS_j$ remains the same, and the process is numerically as accurate as possible.  2)  The cost for computing  the product $S_j^{-1}AS_j$  is only $O(n)$. Thus, a gain of an order-of- magnitude is guaranteed. 3) The product $S_j^{-1}AS_j$  does not fills-up the matrix and preserves all the created zeros in previous steps. Also, to restore  the $J$-Hessenberg form of $A$, only  a cheaper  additional cost of $O(n^2)$ is needed.
 
  In the sequel, we explain first how one may remedy to the fatal breakdown, encountered by JHESS, when applied to the  example $A_6$ and highlights the main lines of the method. Then we present a method to cure the fatal breakdown in  the general case. The idea is the following :  one seeks for a symplectic transforming matrix $S$ so that the similar matrix $S A_6 S^{-1}$, may be reduced by JHESS.  The choice of $S$ should be done carefully. A judicious choice of $S$ consists in taken $S$ equal to Van Loan's Householder  matrix 
  \begin{equation}\label{vanh1}
  S=\left(
  \begin{array}{llll}
  H_2 & 0 & 0&0\\
  0 & 1 & 0 & 0\\
  0 & 0& H_2 & 0\\
  0 & 0 & 0 & 1
  \end{array}
  \right),
  \end{equation}
 or Van Loan's Givens  matrix 
 \begin{equation}\label{dirg1}
 S=\left(
 \begin{array}{llll}
 G_2 & 0 & 0&0\\
 0 & 1 & 0 & 0\\
 0 & 0& G_2 & 0\\
 0 & 0 & 0 & 1
 \end{array}
 \right),
 \end{equation}
  where $H_2$ (respectively $G_2$) is a 2-by-2 Householder matrix (respectively  a 2-by-2 Givens matrix) such that $H_2 (1,2)^T= \sqrt{5}(1,0)^T$  (respectively $G_2 (1,2)^T= \sqrt{5}(1,0)^T$).  If we proceed with choices (\ref{vanh1}) or (\ref{dirg1}), we get the first column of $SA_6$  proportional to $e_1$ and only rows 1,2 and 4,5 of $SA_6$ may change. With the choice (\ref{dirg1}), we obtain $G_2=\left(
  \begin{array}{ll}
  c & s\\
  -s & c
  \end{array}
  \right),$ with $c=1/\sqrt{5},$ $s=2c,$ and 
  $$SA_6 = \left(
  \begin{array}{llllll}
  \sqrt{5} & 2/\sqrt{5} & 0 & \sqrt{5} & 4/\sqrt{5} & 0\\
   0 & 1/\sqrt{5} & 0 & 0& -3/\sqrt{5} & 0\\
    0 & 2 & 1 & 0& 2 & 1\\
     0& 4/\sqrt{5} & 4/\sqrt{5} & 7/\sqrt{5} & 2/\sqrt{5} & 0\\
       0& -3/\sqrt{5} & 2/\sqrt{5} & 1/\sqrt{5} & 1/\sqrt{5} & 0\\
       0 & 0 & 1 & 0 & 3 & 1
  \end{array}
  \right).$$
  The multiplication of  $SA_6$ on the left by $S^{-1}$  acts only on the columns 1, 2 and 4, 5 of 
   $SA_6.$ The other columns remain unchanged. We obtain 
   \begin{equation}\label{exam2}
   SA_6S^{-1}= \left(
   \begin{array}{llllll}
   9/5 & -8/5 & 0 & 13/5 & -6/5 & 0\\
   2/5 & 1/5 & 0 & -6/5 & 3/5 & 0\\
   4/\sqrt{5 }& 2/\sqrt{5 } & 1 &  4/\sqrt{5 } &  2/\sqrt{5 } & 1 \\
   8/5 & 4/5&  4/\sqrt{5 }& 11/5& _12/5&0\\
   -6/5&-3/5&2/\sqrt{5 } & 3/5& -1/5 & 0\\
   0&0&1&6/\sqrt{5}&3/\sqrt{5} & 1
   \end{array}
  \right).
   \end{equation}
   We applied JHESS (also JHMSH) to the matrix $SA_6 S^{-1}$ of (\ref{exam2}).  The algorithm run well   and the reduction to the $J$-Hessenberg form is obtained.  The SR-algorithm is then applied with  explicit and implicit versions, and both converge. Thus the fatal breakdown of JHESS  (or similarly JHMSH) is cured. 
Recall that the algorithm JHMSH as described in \cite{Sal6} is as follows
 \begin{alg}
 	function [S,H]=JHMSH(A)\\
 	\indent $twon=size(A(:,1));\;n=twon/2;\;
 	S=eye(twon);$\\
 	\indent for $ j=1:n-1$\\
 	\indent\indent    $J=[zeros(n-j+1),eye(n-j+1);-eye(n-j+1),zeros(n-j+1)];$\\
 	\indent\indent    $ro=[j:n,n+j:2n];\;co=[j:n,n+j:2n];$\\
 	\indent\indent    $[c,v]=osh2(A(ro,j));$\\
 	\indent \% Updating $A:$\\
 	\indent\indent    $A(ro,co)=A(ro,co)+c * v * (v' * J * A(ro,co));$\\
 	\indent\indent    $A(:,co)=A(:,co)-(A(:,co) * (c * v)) * v' * J;$\\
 	\indent\% Updating $S$ (if needed):\\
 	\indent\indent    $S(:,co) = S(:,co)- c * (v * v') * J * S(:,co);$ 
 	\\
 	\indent \indent for $k=2n:n+j+1,$\\
 	\indent \indent  $[c,s]=vlg(k,A(:,n+j)),\;$\\
 	\indent	\indent \%Updating $A$:\\
 	\indent	\indent $\left[
 	\begin{array}{l}
 	A(k,co) \\
 	A(n+k,co)	
 	\end{array}
 	\right]=
 	\left[
 	\begin{array}{ll}
 	c & s\\
 	-s & c	
 	\end{array}
 	\right]
 	\left[
 	\begin{array}{l}
 	A(k,co) \\
 	A(n+k,co)	
 	\end{array}
 	\right];$\\
 	\indent	\indent $\left[
 	\begin{array}{ll}
 	A(:,k) &	A(:,n+k)	
 	\end{array}
 	\right]=
 	\left[
 	\begin{array}{ll}
 	A(:,k) &	A(:,n+k)	
 	\end{array}
 	\right]  		
 	\left[
 	\begin{array}{ll}
 	c & -s\\
 	s & c	
 	\end{array}
 	\right];
 	$\\
 	\indent \%Updating $S$ (if needed):\\
 	\indent	\indent $\left[
 	\begin{array}{ll}
 	S(:,k) &	S(:,n+k)	
 	\end{array}
 	\right]=
 	\left[
 	\begin{array}{ll}
 	S(:,k) &	S(:,n+k)	
 	\end{array}
 	\right]  		
 	\left[
 	\begin{array}{ll}
 	c & -s\\
 	s & c	
 	\end{array}
 	\right];
 	$\\ 	
 	\indent\indent end\\
 	\indent \indent if $j\leq n-2$\\
 	\indent\indent	$[\beta,w]=vlh(j+1,A(:,n+j));$\\
 	\indent	\indent \%Updating $A$:\\
 	\indent\indent    $A(j+1:n,co)=A(j+1:n,co)-\beta *w*w'*A(j+1:n,co) $\\
 	\indent\indent    $A(j+1+n:2n,co)=A(j+1+n:2n,co)-\beta *w*w'*A(j+1+n:2n,co);$\\	
 	\indent\indent    $A(:,j+1:n)=A(:,j+1:n)-\beta *A(:,j+1:n)w*w';$\\	
 	\indent\indent    $A(:,n+j+1:2n)=A(:,n+j+1:2n)-\beta *A(:,n+j+1:n)w*w';$\\	
 	\indent \%Updating $S$ (if needed):\\
 	\indent\indent    $S(:,j+1:n)=S(:,j+1:n)-\beta *S(:,j+1:n)w*w';$\\	
 	\indent\indent    $S(:,n+j+1:2n)=S(:,n+j+1:2n)-\beta *S(:,n+j+1:n)w*w';$\\	
 	\indent \indent end	\\
 	\indent \indent end\\
 	\indent 	end\\
 \end{alg}
The breakdown in JHMSH occurs exactly in the same conditions as in JHESS, and is located in the call of the function $osh2.$  A slight different version of JHMSH is JHMSH2 (see \cite{Sal6}).
\subsection{Curing breakdowns  in JHESS, JHMSH algorithms}
We present here, in a general manner,  the strategy of curing breakdowns or near breakdowns   which may occur in JHESS or JHMSH algorithms. Let us remark that breakdowns (or near-breakdowns) in JHESS (respectively in JHMSH) may occur only when the function $gal$  (respectively $osh2$) is called, and hence it concerns only columns from the first half of the current matrix.

	Let  $A \in \mathbb{R}^{2n\times 2n}$ be a  matrix  and let $H$ be the matrix that one  obtains  at stage $1 \leq j \leq n-1$ of the algorithm JHESS, by executing instructions {\bfseries{1.}} to {\bfseries{7.}} (corresponding to the current updated  matrix $A$ in the process, at stage $j$ and until instruction {\bfseries{7.}}).  
	Partitioning 
	$H=\left(
	\begin{array}{ll}
	H_{11} & 	H_{12} \\
	H_{21} & 	H_{22} \\
	\end{array}
	\right),$ then $H_{11},\; 	H_{12},\; 	H_{21},\; 	H_{22}$ have the form  \\
	$$H_{11}=
	\left[
	\begin{array}{lllllll}
	h_{1,1}&h_{1,2}& \ldots & h_{1,j} & h_{1,j+1} & \ldots & h_{1,n} \\
	0&h_{2,2}& \ldots & h_{2,j} & h_{2,j+1} & \ldots & h_{2,n} \\
	\vdots &\ddots& \ddots & \vdots & \vdots& \ldots & \vdots \\
	0 && 0 & h_{j,j} & h_{j,j+1} & \ldots & h_{jn} \\
	0 &\ldots& 0 &  h_{j+1,j} & h_{j+1,j+1} & \ldots & h_{j+1,n}     \\
	0 &\ldots& 0 &  0 & h_{j+2,j+1} & \ldots & h_{j+2,n}     \\
	\vdots &&  \vdots& \vdots & \vdots&  & \vdots   \\
	0 &\ldots& 0 & 0 & h_{n,j+1} & \ldots & h_{n,n} 
	\end{array}
	\right],
	$$
	$$H_{12}=
	\left[
	\begin{array}{lllllll}
	h_{1,n+1}&h_{1,n+2}& \ldots & h_{1,n+j-1} & h_{1,n+j} 
	&\ldots & h_{1,2n} \\
	h_{2,n+1}&h_{2,n+2}& \ldots & h_{2,n+j-1} & h_{2,n+j} 
	& \ldots & h_{2,2n} \\
	\vdots &\ddots& \ddots & \vdots & \vdots 
	&  & \vdots \\
	\vdots && \ddots & h_{j-1,n+j-1} & h_{j-1,n+j}  
	& \ldots & h_{j-1,2n} \\
	0 &\ldots& \ldots & h_{j,n+j-1} & h_{j,n+j}  
	& \ldots & h_{j,2n}     \\
	0 &\ldots& \ldots & 0 & h_{j+1,n+j} 
	& \ldots & h_{j+1,2n}     \\
	0 &\ldots& \ldots & 0 & h_{j+2,n+j}  
	& \ldots & h_{j+2,2n}     \\
	\vdots &&  & \vdots & \vdots 
	&  & \vdots   \\
	0 && \ldots & 0 &  h_{n,n+j}  
	& \ldots & h_{n,2n} 
	\end{array}
	\right],
	$$
	$$H_{21}=\left[
	\begin{array}{lllllll}
	 h_{n+1,1}&   h_{n+1,2}& \ldots & h_{n+1,j} & h_{n+1,j+1} & \ldots & h_{n+1,n}\\
	0 &\ddots& \ddots & \vdots & \vdots& \ldots & \vdots \\
	\vdots &\ddots&  \ddots& h_{n+j-1,j} & h_{n+j-1,j+1} & \ldots & h_{n+j-1,n} \\
	\vdots && \ddots &  h_{n+j,j} & h_{n+j,j+1} & \ldots & h_{n+j,n}\\
	0 &\ldots& \ldots & 0 & h_{n+j+1,j+1} & \ldots & h_{n+j+1,n} \\
	\vdots &&  & \vdots & \vdots&  & \vdots\\
	0 &\ldots& \ldots & 0 & h_{2n,j+1} & \ldots & h_{2n,n} 
	\end{array}
	\right],$$
	$$H_{22}=
	\left[
	\begin{array}{lllllll}
	h_{n+1,n+1}&h_{n+1,n+2}& \ldots  & h_{n+1,n+j-1} & h_{n+1,n+j}  
	&\ldots & h_{n+1,2n} \\
	0&h_{n+2,n+2}& \ldots  &h_{n+2,n+j-1} & h_{n+2,n+j} 
	& \ldots & h_{n+2,2n} \\
	0 &0& \ddots  &  &\vdots
	&  & \vdots \\
	\vdots&\vdots& \ddots  &h_{n+j-1,n+j-1} &h_{n+j-1,n+j} 
	& \ldots & h_{n+j-1,2n}     \\
	0&\ldots& \ldots  &0&h_{n+j,n+j} 
	& \ldots & h_{n+j,2n}     \\
	\vdots &&  &\vdots & \vdots 
	 & & \vdots   \\
	0 &\ldots& \ldots &0 & h_{2n,n+j} 
	& \ldots & h_{2n,2n} 
	\end{array}
	\right].
	$$
The breakdown occurs in JHESS when the coefficient $h_{n+j,j}=0$ and $h_{j+1,j} \neq 0.$ In this case, JHESS stops computations. To overcome this fatal breakdown, we construct the orthogonal matrix 
$P^{(j)}=diag(I_{j-1},H_2^{(j)},I_{n-j-1})$ where $H_2^{(j)}$ is a 2-by-2 Householder matrix. We set $S^{(j)}= diag(P^{(j)}, P^{(j)}).$ The matrix $S^{(j)}$ is symplectic and orthogonal. The choice of the 2-by-2 Householder matrix  $H_2^{(j)}$ is so that $H_2^{(j)}\left( \begin{array}{l}
h_{j,j}\\
h_{j+1,j}
\end{array}
\right)
=
\left(  
\begin{array}{l}
h'_{j,j}\\
0
\end{array}
\right)
.$
Thus, the action $S^{(j)}H$ annihilates the position $(j+1,j)$ of the updated matrix $H$  and keep unchanged all zeros created previously except potentially the position $(j+1,n+j-1)$ (in the block $H_{12})$). Keep in mind that the action of $S^{(j)}H$  on $H$ affects only rows $j,j+1,n+j,n+j+1.$ The action $S^{(j)}H[S^{(j)}]^{-1}$ on $S^{(j)}H$ affects only the columns $j,j+1,n+j,n+j+1$.  Thus, the columns $1,\ldots,j-1$ and $n+1,\ldots,n+j-1$ remain unchanged (hence all zeros created previously are not affected). The process is pursued as follows : one must annihilate the potentially nonzero   entry in position  $(j+1,n+j-1)$  and keep unchanged all zeros created previously in columns $1,\ldots,j-1$ and $n+1,\ldots,n+j-1.$ This may be addressed  by applying the similarity $T^{(j)}H[T^{(j)}] ^{-1}$ to the obtained matrix $H$, where $T^{(j)}=
 diag(Q^{(j)}, Q^{(j)})$
 states for Van Loan's Householder matrix, given by 
$Q^{(j)}=diag(I_{j-1},K_2^{(j)},I_{n-j-1})$ where $K_2^{(j)}$ is a 2-by-2 Householder matrix. The choice of the 2-by-2 Householder matrix  $K_2^{(j)}$ is so that $K_2^{(j)}\left( \begin{array}{l}
h_{j,n+j-1}\\
h_{j+1,n+j-1}
\end{array}
\right)
=
\left(  
\begin{array}{l}
h'_{j,n+j-1}\\
0
\end{array}
\right)
.$
The cost of this curing strategy step is $O(n).$
  Next, the algorithm is pursued normally by executing again step $j.$
\begin{rmk} 
			The 2-by-2 Householder matrix $H_2^{(j)}$ (respectively $K_2^{(j)}$ ) may be replaced by a 2-by-2 Givens matrix $G_2^{(j)}=\left(
	\begin{array}{ll}
	c & s\\
	-s & c
	\end{array}	
		\right)$ 
		 (respectively by a 2-by-2 Givens $L_2^{(j)}=\left(
		 \begin{array}{ll}
		 c' & s'\\
		 -s' & c'
		 \end{array}	
		 \right)$ )
		where the coefficient $c,s$ (respectively $c',s'$) are chosen so that $G_2^{(j)}\left( \begin{array}{l}
		h_{j,j}\\
		h_{j+1,j}
		\end{array}
		\right)
		=
		\left(  
		\begin{array}{l}
		h'_{j,j}\\
		0
		\end{array}
		\right)
		$  ( respectively $L_2^{(j)}\left( \begin{array}{l}
		h_{j,n+j-1}\\
		h_{j+1,n+j-1}
		\end{array}
		\right)
		=
		\left(  
		\begin{array}{l}
		h'_{j,n+j-1}\\
		0
		\end{array}
		\right)
		$).
		\end{rmk}
		 It is worth noting that one may take an arbitrary  $k$-by-$k$ Householder matrix $H_k^{(j)}$ instead of $H_2^{(j)}$, with  
		 $P^{(j)}=diag(I_{j-1},H_k^{(j)},I_{n-j-k+1}).$ The action $P^{(j)}H$ on $H$ affects only rows $j,\ldots,j+k-1,$  rows $n+j,\ldots,n+j+k-1.$  The action $P^{(j)}H[P^{(j)}]^{-1}$ on $P^{(j)}H$ affects only columns $j,\ldots,j+k-1,$  and columns $n+j,\ldots,n+j+k-1.$ Hence, all zeros created previously in columns $1,\ldots,j-1$ and columns $n+1,\ldots, n+j-1$ remain unchanged, except potentially in positions $(j+1,n+j-1) ,\ldots, (j+k-1,n+j-1).$
		  To pursue the process,  one must annihilate the potentially nonzero   entries in position  $(j+1,n+j-1),\ldots, (j+k-1,n+j-1)$  and keep unchanged all zeros created previously in columns $1,\ldots,j-1$ and $n+1,\ldots,n+j-1.$ This may be addressed  by applying the similarity $T^{(j)}H[T^{(j)}] ^{-1}$ to the obtained matrix $H$, where $T^{(j)}=
		  diag(Q^{(j)}, Q^{(j)})$
		  states for Van Loan's Householder matrix, given by 
		  $Q^{(j)}=diag(I_{j-1},K_k^{(j)},I_{n-j-1})$ where $K_k^{(j)}$ is a $k$-by-$k$ Householder matrix. The choice of the $k$-by-$k$ Householder matrix  $K_k^{(j)}$ is so that $K_k^{(j)}\left( \begin{array}{l}
		  h_{j,n+j-1}\\
		  \vdots\\
		  h_{j+k-1,n+j-1}
		  \end{array}
		  \right)
		  =
		  \left(  
		  \begin{array}{l}
		  h'_{j,n+j-1}\\
		  0
		  \end{array}
		  \right)
		  .$
		  The cost of this curing strategy step is $O(kn).$
		  \subsection{Curing  near-breakdowns in JHESS, JHMSH algorithms}
The near-breakdown occurs in JHESS (or in JHMSH)  when the coefficients $h_{n+j,j}$ and $h_{j+1,j} $  
 are both different from zero  but are near  to the situation of breakdown. This can be measured by the fact that the ratio $ \displaystyle{
\frac{h_{j+1,j}}
{h_{n+j,j}}}$ is very large. 
 In this case, the non-orthogonal and symplectic transformations involved in JHESS (respectively JHMSH) become ill-conditioned and numerical instability is encountered reducing the accuracy of the reduction. 	In order to remedy to a such near breakdown in JHESS (or JHMSH) algorithm, one may proceed exactly as for curing  a breakdown, the only difference is that the test $h_{n+j,j}=0$ and $h_{j+1,j} \neq 0$   (corresponding to a breakdown) is replaced by the $ \displaystyle{
 	\frac{h_{j+1,j}}
 	{h_{n+j,j}}} \geq \tau$ (corresponding to a near-breakdown), where $\tau$ is a certain tolerance.
 \subsection{SR algorithm}
 The SR algorithm, is a QR like algorithm which can roughly be described as follows. For a given matrix $M \in \mathbb{R}^{2n\times 2n}$, it computes:\\
 1. An upper $J$-Hessenberg reduction $M_1=S_0^J M S_0$ where $S_0$ is symplectic. Set $S=S_0.$\\
 2. Iteration : For $k=1,\ldots,$ compute $M_{k+1}=S_k^{J} M_k S_{k}$ where $S_k$ stands for the symplectic factor of the SR decomposition $p_k(M_k) = S_k R_k$ of a polynomial $p_k$ of $M_k$ and update $S= S S_k.$ \\
 The iterate $M_{k+1}$ remain $J$-Hessenberg if the  matrix $M_k$ is $J$-Hessenberg. Like the QR algorithm, SR algorithm admits an implicit version : the decompositions $p_k(M_k) = S_k R_k$ are not performed explicitly. Since SR algorithm is based on  $J$-Hessenberg reductions and SR decompositions, breakdowns or near-breakdowns may be encountered both in the explicit or implicit versions of the algorithm. Of course, the implicit form is preferred to the explicit one.  In \cite{Buns} there is no strategy proposed when a breakdown is meet. The algorithm is topped. However, a  technique has been  proposed  in the situation of a near breakdown, occurring at the iteration $j.$ 
 %
 
 The proposed technique consists in computing the similarity $M_{j+1}=S_j M_j[S_j]^{-1}$, where 
 $S_j = I -ww^T J$ and 
 $w\in 
 \mathbb{R}^{2n\times 2n}$ is a random vector, with $\norm{w}_2=1.$ The algorithm continue with $M_{j+1}.$ This technique presents several serious drawbacks : \\
 1. The symplectic transformation $(I-ww^T J)$ is never  orthogonal (except for $w=0$) and hence  its condition number may be large. The condition number of $M_{j+1}$ could be worse than this of  $M_{j}.$  \\
 2. Computing the similarity  $M_{j+1}=S_j M_j[S_j]^{-1}$ costs $O(n^2).$\\
 3. The similarity $M_{j+1}=S_j M_j[S_j]^{-1}$  destroys the structure $J$-Hessenberg of $M_j.$ Thus $M_j$ is no longer $J$-Hessenberg. Moreover, the matrix $M_{j+1}$ fills up. Hence, a reduction to a $J$-Hessenberg form is needed to restore the previous structure. The cost of this restoration is $O(n^3)$ which is very expensive. Instead of this similarity, when  breakdown or near breakdown occurs with respect to the column say $i$ of the matrix $M_j$, we propose the similarity 
  $M_{j+1}=P^{(j)} M_j[P^{(j)}]^{-1},$  where   $P^{(j)}=diag(I_{i-1},H_l^{(i)},I_{n-i-l+1})$  and 
  $H_l^{(i)}$ an arbitrary  $l$-by-$l$ Householder matrix. 
   The action $P^{(j)}M_j$ on $M_j$ affects only rows $i,\ldots,i+l-1,$  rows $n+i,\ldots,n+i+l-1.$  The action $P^{(j)}M_j[P^{(j)}]^{-1}$ on $P^{(j)}M_j$ affects only columns $i,\ldots,i+l-1,$  and columns $n+i,\ldots,n+i+l-1.$ Hence, all zeros in columns $1,\ldots,i-1$ and columns $n+1,\ldots, n+i-1$ (because of the form $J$-Hessenberg of $M_j$) remain unchanged, except potentially in positions $(i+1,n+i-1) ,\ldots, (i+l-1,n+i-1).$ To pursue the process, one restores the $J$-Hessenberg form. 
    The advantage of this similarity is that first $P^{(j)}$ is symplectic and orthogonal (hence it is stable), preserves most of the created zeros because of the $J$-Hessenberg structure of $M_j$ and the cost of restoring the $J$-Hessenberg form of $M_j$ do not exceed $O(i n).$
 Unlike QR algorithm, SR algorithm still needs profound investigations. This will be the aim of a forthcoming work.
 \section{Numerical experiments}
 To illustrate our purpose, we consider the following numerical example. Let $A$ be the 12-by-12 matrix
 \[
 A =
 \begin{pmatrix}
 \begin{array}{rrrrrrrrrrrr}
 1 & 5 & 7 & 9 & 5 & 1 & 1 & 3 & 1 & 3 & 7 & 2 \\
 0 & 1 & 4 & 6 & 1 & 2 & 2 & 1 & 5 & 4 & 3 & 5 \\
 0 & 0 & 1 & 2 & 3 & 2 & 0 & 0 & 1 & 2 & 5 & 3 \\
 0 & 0 & 2 & 1 & 9 & 8 & 0 & 0 & 2 & 1 & 2 & 4 \\
 0 & 0 & 0 & 2 & 1 & 3 & 0 & 0 & 5 & 2 & 1 & 2 \\
 0 & 0 & 0 & 4 & 2 & 1 & 0 & 0 & 4 & 3 & 2 & 1 \\
 1 & 4 & 7 & 2 & 1 & 3 & 1 & 7 & 6 & 1 & 6 & 7 \\
 0 & 1 & 9 & 3 & 5 & 1 & 0 & 1 & 4 & 5 & 8 & 3 \\
 0 & 0 & 0 & 2 & 7 & 9 & 0 & 0 & 1 & 3 & 4 & 5 \\
 0 & 0 & 0 & 1 & 2 & 8 & 0 & 0 & 3 & 1 & 7 & 3 \\
 0 & 0 & 0 & 2 & 1 & 2 & 0 & 0 & 4 & 3 & 1 & 2 \\
 0 & 0 & 0 & 9 & 3 & 1 & 0 & 0 & 1 & 2 & 3 & 1
 \end{array}
 \end{pmatrix}. 
 \]
Following the steps of the algorithms JHESS, JHMSH and JHMSH2, one remarks that
 the condition of a breakdown is fulfilled at the beginning of the step $j=3$ for all of them.
  Let us call MJHESS (respectively JHM$^2$SH and JHM$^2$SH2) the modified algorithm JHESS (respectively JHMSH and JHMSH2) obtained by applying our strategy for curing breakdowns. We obtain the following numerical results, showing the efficiency of the method. 
   \begin{table*}[h]
   	\centering
   	\begin{tabular}{lllll}
   		\toprule
   		\multirow{2}*{$2n$} & \multicolumn{4}{c}{Loss of $J$-Orthogonality $\left\Vert I-S^{J}S\right\Vert _{2}$}\\
   		\cmidrule{2-5} 
   		& \multicolumn{1}{c}{$ JHESS $} & \multicolumn{1}{c}{$ MJHESS $} & \multicolumn{1}{c}{$ JHM^{2}SH $} 
   		& \multicolumn{1}{c}{$ JHM^{2}SH2 $} \\ 
   		\midrule
   		$12$ & fails & $1.8553e-15$ & $5.0842e-15$ & $6.6428e-15$ \\ 
   		\bottomrule
   	\end{tabular}
   \end{table*}
  Thus the $J$-orthogonality is numerically preserved up to the machine precision for MJHESS (respectively JHM$^2$SH and JHM$^2$SH2). It is worth noting that preserving the $J$-orthogonality is  crucial  for SR-algorithm in order to get accurate eigenvalues and vectors of a matrix. 
   \begin{table*}[h]
   	\centering
   	\begin{tabular}{lllll}
   		\toprule
   		\multirow{2}*{$2n$} & 
   		\multicolumn{4}{c}{Reduction error  $ \left\Vert A-SHS^{-1} \right\Vert _{2} $}\\
   		\cmidrule{2-5}
   		& \multicolumn{1}{c}{$ JHESS $} & \multicolumn{1}{c}{$ MJHESS $} & \multicolumn{1}{c}{$ JHM^{2}SH $} 
   		& \multicolumn{1}{c}{$ JHM^{2}SH2 $} \\
   		\midrule
   		$12$ & fails & $3.2709e-14$ & $3.8777e-13$ & $2.7653e-13$ \\ 
   		\bottomrule
   	\end{tabular}
   \end{table*}
  
   One observes also that the error in the reduction to $J$-Hessenberg form is very satisfactory for MJHESS (respectively JHM$^2$SH and JHM$^2$SH2). Notice that the algorithm JHESS  as given in \cite{Buns}, applied to the matrix $A$,  without our strategy for curing breakdown, simply fails to perform a reduction to a $J$-Hessenberg form. 

 Let $A=S_i H_i S_i^{-1}, \;i=1,2,3$ be the $J$-Hessenberg reduction  obtained respectively by the algorithms MJHESS, JHM$^2$SH and JHM$^2$SH2. It is known (see \cite{Buns}) that there exist $D_1,\;D_2,\;D_3$ such that 
 $S_1 D_1 = S_2,\;D_1^{-1}H_1 D_1=H_2,\;S_1 D_2 = S_3,\;D_2^{-1}H_1 D_2=H_3,\;S_2 D_3=S_3,$ and $D_3^{-1}H_2 D_3=H_3,$ with each matrix 
 $D_{i}=\begin{pmatrix}
 \begin{array}{ll}
 C_i & F_i\\
 0 & C_i^{-1}
 \end{array}
 \end{pmatrix},$ where $C_i$ and $F_i$ are diagonals. 
 
We obtain numerically  
$C_1=diag(1,1,0.4113,0.3688,0.7747,0.6638)$ and $F_1=diag(0,0,-2.3962,-2.3371,1.2880,-1.9982),$
and 
$$ \Vert D_{1}^{-1} H_{1} D_{1} - H_{2} \Vert _{2} =  = 5.3417e-13.$$
Also, we have\\ $C_2=(1,1,0.4113,0.3688,0.7747,0.6638),\; F_2=diag(0,0,-2.3962,-2.3371,1.2880,-1.9982),$\\
and 
$$ \Vert D_{2}^{-1} H_{1} D_{2} - H_{3} \Vert _{2} = 3.7881e-13,$$
and finally 
\[
D_{3} = 
I_{12},
\]  where $I_{12}$ stands for identity matrix of size 12, 
with
$$ \Vert D_{3}^{-1} H_{2} D_{3} - H_{3} \Vert _{2}  = 9.4799e-13.$$
Thus the  matrices $D_i$ have numerically the desired forms and as expected, the algorithms  JHM$^2$SH and JHM$^2$SH2 perform quite the same results.
\section{Conclusions}
In this work, we linked the necessary and sufficient condition of the existence of a SR-decomposition with the computations during the process, of  some coefficients of the current matrix. The SR-decomposition is intimately related to the $J$-Hessenberg reduction via the algorithm JHESS. The later (also JHMSH and its different variants) may encounter  fatal breakdowns or suffer from  near-breakdowns. We derive efficient strategies for treating them. The numerical experiments show the efficiency of these strategies.


\begin{thebibliography}{00}

\bibitem{Art}%
E.~Artin, \emph{Geometric Algebra}, Interscience Publishers,
 New York, 1957.




\bibitem{Buns}
A.~Bunse-Gerstner and V.~Mehrmann, {A symplectic {QR}-like
algorithm for
  the solution of the real algebraic {Riccati} equation}, IEEE Trans. Automat.
  Control \textbf{AC-31} (1986), 1104--1113.

 \bibitem{Buns1}%
A.~Bunse-Gerstner, {Matrix factorizations for symplectic {QR}-like
methods}, Linear
  Algebra Appl. \textbf{83} (1986), 49--77.

\bibitem{Dela2}%
J.~Della-Dora, {Numerical linear algorithms and group theory},
Linear
  Algebra Appl. \textbf{10} (1975), 267--283.


\bibitem{Els}
L. ~Elsner, {On some algebraic problems in connection with general
eigenvalue algorithms}, Linear Alg. Appl., 26:123-38 (1979).

\bibitem{Golu}
G.~Golub and C.~Van Loan, \emph{Matrix {Computations}}, third ed.,
The Johns
  Hopkins U.P., Baltimore, 1996.








\bibitem{Paig}%
C.~Paige and C.~Van Loan, {A {Schur} decomposition for
{Hamiltonian}
  matrices}, Linear Algebra Appl. \textbf{41} (1981), 11--32.


  \bibitem{Sal1}%
  A.~Salam, {On  theoretical and numerical aspects of symplectic Gram-Schmidt-like
  algorithms}, Numer. Algo., \textbf{39} (2005), 237-242.



\bibitem{Sal2}
A. Salam, A. El Farouk, E. Al-Aidarous, Symplectic Householder Transformations for a QR-like decomposition, a Geometric and Algebraic Approaches, J. of Comput. and Appl. Math., Vol. 214, Issue 2, 1 May 2008, Pages 533-548.
\bibitem{Sal3}
A. Salam and E. Al-Aidarous and A. Elfarouk, Optimal symplectic Householder transformations for SR-decomposition, Linear Algebra and Its Appl., 429 (2008), no. 5-6, 1334-1353.
\bibitem{Sal4}
A. Salam, E. Al-Aidarous,  Error analysis and computational aspects of SR factorization, via optimal symplectic Householder Transformations, Electronic Trans. on Numer. Anal., Vol. 33, pp. 189-206, 2009.
\bibitem{Sal6}
A. Salam and H. Ben Kahla, An   upper $J$- Hessenberg reduction of a matrix through symplectic Householder  transformations, submitted.

\bibitem{Vloa}%
 C.~Van Loan, {A symplectic method for approximating all the
eigenvalues of
  a {Hamiltonian} matrix}, Linear Algebra Appl. \textbf{61} (1984), 233--251.

\bibitem{Wat}
D.S.~Watkins,{\it The Matrix Eigenvalue Problem : GR and Krylov subspace methods}, SIAM, 2007.

\bibitem{Wilk}%
J.H.~Wilkinson, \emph{ The {Algebraic} {Eigenvalue} {Problem}},
Clarendon Press, Oxford, England.

\end{thebibliography}
\end{document}